\documentclass[a4paper,11pt,reqno]{amsart}

\usepackage{amssymb}
\usepackage[all]{xy}

\usepackage{amsmath}
\usepackage{amsthm}    
\usepackage{mathrsfs}

\usepackage{hyperref}
\usepackage{xcolor}
\usepackage{enumitem}
\usepackage{subcaption}
\usepackage[final]{graphicx}
\usepackage{float}
\usepackage{epic}
\usepackage{setspace}

\usepackage[all]{xy}
\usepackage{color}
\usepackage{calligra}
\usepackage[T1]{fontenc}

\usepackage{verbatim}

\newcount\mycount

\numberwithin{equation}{section}

\numberwithin{equation}{section}
\newtheorem{introtheorem}{Theorem}

\newtheorem{introcorollary}[introtheorem]{Corollary}

\newtheorem{theorem}{Theorem}[section]

\newtheorem{lemma}[theorem]{Lemma}

\newtheorem{proposition}[theorem]{Proposition}

\newtheorem{corollary}[theorem]{Corollary}
\newtheorem*{theorem*}{Theorem}

\theoremstyle{definition}
\newtheorem{definition}[theorem]{Definition}
\newtheorem*{definition*}{Definition}

\newtheorem{remark}[theorem]{Remark}

\theoremstyle{remark}
\newtheorem{example}[theorem]{Example}

\theoremstyle{remark}

\newcommand{\Cs}{\ensuremath{\mathrm{C}^\ast}}
\newcommand{\Csr}{\mathrm{C}^\ast_{\rm r}}

\newcommand{\Prob}{\ensuremath{\mathrm{Prob}}}

\begin{document}

\title[Pureness and stable rank one for reduced twisted group \Cs-algebras]{Pureness and stable rank one for reduced twisted group \Cs-algebras of certain group extensions}

\author{Felipe Flores}
\author{Mario Klisse}
\author{M\'iche\'al \'O Cobhthaigh}
\author{Matteo Pagliero}

\address{Department of Mathematics, University of Virginia, Kerchof Hall. 141 Cabell Dr, Charlottesville, Virginia, United States}

\email{hmy3tf@virginia.edu}

\email{des5sf@virginia.edu}

\address{Department of Mathematics, Christian-Albrechts University Kiel,\\ Heinrich-Hecht-Platz 6,
Kiel, Germany}

\email{klisse@math.uni-kiel.de}

\address{Department of Mathematics and Computer Science, University of Southern Denmark, Campusvej 55,
Odense, Denmark}

\email{pagliero@imada.sdu.dk}

\date{\today. }

\begin{abstract}
We present a sufficient condition for a (twisted) group action on a $\mathrm{C}^\ast$-algebra such that the corresponding  reduced (twisted) group $\mathrm{C}^\ast$-algebra inclusion into the  reduced (twisted) crossed product is selfless. This condition is an action-dependent version of Ozawa's $\mathrm{PHP}$ property for groups. We show that topologically free extreme boundary actions and projective actions of lattices in ${\rm PSL}(n,\mathbb R)$ have this property, and thus we provide examples of selfless inclusions from crossed products. As a byproduct, we also obtain that reduced (twisted) group $\mathrm{C}^\ast$-algebras of some group extensions of the form finite-by-$G$, with $G$ having the property $\mathrm{PHP}$, have stable rank one and are pure, and in particular have strict comparison. Examples include all acylindrically hyperbolic groups and all lattices in ${\rm SL}(n,\mathbb R)$ for $n\geq2$.
\end{abstract}

\maketitle


\section*{Introduction}


The notion of {\it strict comparison of positive elements}, or {\it strict comparison} for short, was introduced by Blackadar in \cite{Bl88} as a \Cs-analogue of the property of II$_1$-factors that Murray--von Neumann comparison of projections is determined by their trace. We refer the reader to \cite[Definition 3.3]{APTV24} for the definition of strict comparison, and to the same article for a detailed discussion on this topic.
Suffice it to say that it is a regularity property for \Cs-algebras that plays a prominent role in both the Toms--Winter conjecture and Elliott's classification program (see, for example, \cite{MS12,Wi18,Th20}), and has also found applications in continuous model theory \cite{KES25}, dynamical systems \cite{Ke20}, and time-frequency analysis \cite{BEvV22}.

Establishing strict comparison is usually a difficult problem, especially for nonnuclear C$^\ast$-algebras. This is well exemplified by the fact that, until recently, it was an open problem whether the reduced group \Cs-algebra of the free group $\mathbb{F}_n$ has strict comparison when $2 \leq n < \infty$. 
The situation drastically changed thanks to the breakthrough article of Amrutam, Gao, Kunnawalkam Elayavalli, and Patchell \cite{AGKEP25}, in which they proved that reduced group \Cs-algebras of finitely generated acylindrically hyperbolic groups with trivial finite radical and satisfying the rapid decay property are \emph{selfless} in the sense of Robert when equipped with the canonical trace. 

\begin{definition*}[see {\cite[Definition 2.1 and Theorem 2.6]{Ro25}, \cite[Section 6]{Oz25}, \cite[Definition 2.3]{HKEPR25}}]
Let $(A,\rho)$ be a \Cs-probability space, i.e., a unital \Cs-algebra equipped with a GNS-faithful state. Then $(A,\rho)$ is said to be \emph{selfless}
if $A\not\cong\mathbb{C}$ and there exists a unital \Cs-algebra $B\not\cong\mathbb C$ equipped with a GNS-faithful state $\psi$ such that the first factor embedding into the reduced free product \Cs-algebra $\iota\colon (A,\rho) \hookrightarrow (A,\rho) \star (B,\psi)$ is \emph{existential}, that is, there exists an ultrafilter $\mathcal U$ and a state-preserving embedding $\theta\colon(A,\rho)\star(B,\psi)\hookrightarrow(A^{\mathcal U},\rho^{\mathcal U})$ such that $\theta \circ \iota$ agrees with the diagonal embedding.

Moreover, $(A,\rho)$ is said to be \emph{completely selfless} if for every \Cs-algebra $C$, $\theta$ induces a $\ast$-homomorphism 
\[
((A,\rho)\star(B,\psi))\otimes C \longrightarrow ((A,\rho)\otimes C)^{\mathcal U} .
\]

Also, if $A_0\subseteq (A,\rho)$ is an inclusion of \Cs-probability spaces, i.e., $A_0$ is a \Cs-subalgebra of $A$ that when equipped with the restriction of $\rho$ is a \Cs-probability space, and $\theta(B)\subseteq A_0^{\mathcal U}$, one says that  $A_0 \subseteq (A,\rho)$ is a \emph{selfless inclusion}. In this case, every intermediate \Cs-probability space is selfless as well. 
\end{definition*}

When equipped with a tracial state, \cite[Theorem 3.1]{Ro25} (see also \cite{BarlakSzabo}) shows that a selfless \Cs-algebra is simple, monotracial with unique quasitracial state, has strict comparison, and stable rank one in the sense of Rieffel \cite{Rie83}, that is, the collection of all invertible elements is dense.

A \Cs-algebra is {\it pure}, according to Winter's original definition (see \cite[Definition 3.6(i)]{Wi12}), if it has strict comparison and it is almost divisible.
Almost divisibility, like strict comparison, is automatic for II$_1$-factors. In a II$_1$-factor, it is the property that for every projection $p$ and $n\geq2$ there exists a projection $q$ on which the unique trace $\tau$ has value $\frac{1}{n}\tau(p)$.
We refer the reader to \cite{Wi12,APTV24,PTV25}  for a comprehensive study of pure \Cs-algebras.
The key result relating pure and selfless \Cs-algebras is \cite[Theorem 5.13]{APTV24}: Every simple, unital, non-elementary \Cs-algebra with a unique normalised quasitrace and strict comparison is pure.
In particular, (completely) selfless tracial \Cs-algebras are pure.
As a result, the main theorem of \cite{AGKEP25} shows that many group \Cs-algebras, including $\Cs_r(\mathbb F_n)$ for $2 \leq n <\infty$, are pure and have stable rank one.
This result has inspired new works on selfless \Cs-algebras arising from constructions such as (twisted) group algebras \cite{Vi25,KEPT25,Oz25,MWY25,RTV25}, graph products and free products \cite{HKER25,HKEPR25,FKOCP25}, and crossed products \cite{GJKEPR26,Ohs26,Bell}.

The first objective of this article is to find new examples of selfless reduced twisted group \Cs-algebras, and give the first instances of selfless inclusions of the reduced  reduced (twisted)  group \Cs-algebra inside a (twisted) crossed product.
Recall that Ozawa showed in \cite{Oz25} that reduced group \Cs-algebras of discrete groups satisfying the $\mathrm{PHP}$ property\footnote{To be precise, Ozawa uses $\mathrm{P}_{\mathrm{PHP}}$ instead of $\mathrm{PHP}$. The name is a tribute to the works of Powers, Haagerup, and Pisier on similar notions.} are completely selfless. It is not hard to see that his proof remains valid in presence of a scalar-valued $2$-cocycle, and thus transfers elegantly to the twisted setting as well. Moreover, the $\mathrm{PHP}$ property is also easily adapted to an inclusion of groups, in which case the respective inclusion of reduced (twisted) group \Cs-algebras is selfless.
In the present work, we introduce an action-dependent adaptation of the $\mathrm{PHP}$ property that we use to upgrade Ozawa's result to a statement about inclusions. To be more precise, by requiring that the action of a $\mathrm{PHP}$ group satisfies an averaging-type condition linked to data provided by the regular $\mathrm{PHP}$ property, we can show that the inclusion of the  reduced (twisted) group \Cs-algebra inside the reduced (twisted) crossed product is selfless. More precisely, we prove the following result. (For the precise definitions of the $\mathrm{PHP}$ property for an action and for an inclusion, see Definition~\ref{def:PHP}.) 

\begin{introtheorem}[{see Theorem~\ref{thm:crossed-product-inclusion}}]\label{MainTh2}
Let $G$ be a discrete group, $A$ a unital \Cs-algebra with a state $\phi$,  $(\alpha,\frak{u})\colon G\curvearrowright A$ a twisted action with the $\mathrm{PHP}$ property with respect to $\phi$, and assume that $\phi \circ E$ and its restriction to $\Cs_r(G,\frak{u})$ are GNS-faithful. Then the inclusion 
\[
\Cs_r(G,\frak{u}) \subseteq (A \rtimes_{r,\alpha,\frak{u}} G,\phi \circ E)
\]
is selfless, and every intermediate \Cs-probability space is completely selfless.

Moreover, if $(\alpha,\frak{u})$ has the  $\mathrm{PHP}$ relative to a subgroup $H\le G$, and the restriction of $\phi \circ E$ to $\Cs_r(H,\frak{u}|_{H})$ is GNS-faithful, the same conclusion holds for the inclusion
\[
\Cs_r(H,\frak{u}|_{H}) \subseteq (A \rtimes_{r,\alpha,\frak{u}} G,\phi \circ E).
\]
\end{introtheorem}

The class of groups with property $\mathrm{PHP}$ includes all groups $G$ with a topologically free extremely proximal action $G\curvearrowright X$ as shown in \cite{Oz25}. By the main theorem of \cite{Yan25}, this includes all acylindrically hyperbolic groups with trivial finite radical. Other examples can be found in \cite{BIO20,IvOm17}.

In Theorem~\ref{thm:extremely-proximal-PHP} we show that the action $G\curvearrowright C(X)$ has, in fact, the $\mathrm{PHP}$ property. 
Moreover, for any finite-index subgroup $H\leq G$, the restricted action $H\curvearrowright X$ has the property $\mathrm{PHP}$ as well (see Proposition \ref{pass}). Hence, we obtain the following corollary.

\begin{introcorollary}[{see~Corollary~\ref{cor:extreme}}]\label{introcor:B}
Let $X$ be a compact Hausdorff space with $|X|>2$.
For any twisted action $(\alpha,\frak{u})\colon G\curvearrowright C(X)$ such that $\alpha$ induces an extremely proximal, topologically free action on $X$, and any $\phi\in S(C(X))$, the inclusion 
\[
\Cs_r(G,\frak{u}) \subseteq (C(X)\rtimes_{r,\alpha,\frak{u}}G,\phi\circ E)
\]
is selfless, and every intermediate \Cs-probability space is completely selfless.

Moreover, if $H\le G$ is a finite index subgroup and $\frak{v}=\frak{u}|_{H}$ the restricted cocycle, the same conclusion holds for the inclusion
\[
\Cs_r(H,\frak{v}) \subseteq (C(X)\rtimes_{r,\alpha,\frak{u}}G,\phi\circ E).
\]
\end{introcorollary}

If a group is selfless in the sense of \cite{AGKEP25} and has rapid decay, the twisted group \Cs-algebra by a scalar-valued $2$-cocycle is selfless by \cite[Theorem A]{RTV25}. This is the case for all acylindrically hyperbolic groups with trivial finite radical and with rapid decay, as they admit an extreme boundary by \cite{Yan25}.
We observe that as a byproduct of Corollary \ref{introcor:B}, the same conclusion is true even when one does \emph{not} assume the rapid decay property, and for any $C(X)$-valued $2$-cocycle. 

Another class of $\mathrm{PHP}$ groups comes from Zariski-dense subgroups $\Gamma$ of ${\rm PSL}(n,\mathbb R)$, for $n\geq 2$; see \cite[Proposition 16]{Oz25}.
In Theorem~\ref{thm:PSL-PHP} we also strengthen this result by establishing that the canonical action $\Gamma\curvearrowright\mathbb P^{n-1}(\mathbb R)$ has the ${\rm PHP}$ property. We thus obtain the following corollary.

\begin{introcorollary}[{see~Corollary~\ref{cor:PSL}}]
Let $n\geq2$, and $\Gamma \le \mathrm{PSL}(n,\mathbb R)$ be a Zariski-dense subgroup such that the canonical action $\Gamma\curvearrowright\mathbb P^{n-1}(\mathbb R)$ is minimal.
Let $\frak{u}$ be any $2$-cocycle for $\alpha$, and $\phi\in S(\mathbb P^{n-1}(\mathbb R))$ any state. Then, the inclusion
\[
\Cs_r(\Gamma,\frak{u}) \subseteq (C(\mathbb P^{n-1}(\mathbb R))\rtimes_{\alpha,\frak{u},r}\Gamma,\phi \circ E)
\]
is selfless, and every intermediate \Cs-probability space is completely selfless.

In particular, the assumptions are automatic whenever $\Gamma\le \mathrm{PSL}(n,\mathbb R)$ is a lattice.
\end{introcorollary}

We also give a sufficient condition to obtain $\mathrm{PHP}$ actions when $G$ acts non-faithfully on a unital \Cs-algebra $A$; see Theorem~\ref{thm:F2-actions}. 
This new criterion yields the following result as a byproduct.

\begin{introcorollary}[{see~Corollary~\ref{cor:strongly-proximal}}]
    Let $G$ be an acylindrically hyperbolic group whose finite radical is trivial, and $X$ a compact Hausdorff space. Let $(\alpha,\frak u)\colon G\curvearrowright C(X)$ be a twisted action whose induced action on $X$ is strongly proximal, minimal, and such that $K:=\ker(\alpha)$ is nontrivial.
    Then, for any state $\phi\in S(C(X))$, the inclusion 
    \[
    \Cs_r( G,\frak{u}) \subseteq (C(X)\rtimes_{\alpha,\frak{u},r} G, \phi \circ E)
    \]
    is selfless, and each intermediate \Cs-probability space is completely selfless. 
\end{introcorollary}

We also mention here that one can proceed exactly as in the proof of \cite[Corollary 2.1]{RTV25}, which uses the structural results from \cite{APTV24,APT18}, and describe the Cuntz semigroup of the \Cs-algebras covered by Theorem \ref{MainTh2} as follows: 

\begin{introcorollary}
    Let $G$ be a countable group with property $\mathrm{PHP}$, and  $\frak{u}\in Z^2(G,\mathbb T)$  a $2$-cocycle. 
    Then the Cuntz semigroup of $\Cs_r(G,\frak u)$ is
    \[
    \mathrm{Cu}(\Cs_r(G,\frak u)) \cong V(\Cs_r(G,\frak u)) \sqcup (0,\infty] \,, 
    \]
    where $V(\Cs_r(G,\frak u))$ denotes the Murray--von Neumann semigroup of $\Cs_r(G,\frak u)$.
\end{introcorollary}

As a second goal of the present article, we continue the line of work started by Raum--Thiel--Vilalta \cite{RTV25}, whose objective is to first transfer the regularity properties mentioned before from the untwisted case to the twisted case, and then to use these results to obtain regularity properties for \Cs-algebras of certain group extensions.
In particular, we consider the reduced (twisted) group \Cs-algebra of a group $\widetilde{G}$ arising from a short exact sequence
$$
1\longrightarrow K\longrightarrow \widetilde{G}\longrightarrow G \longrightarrow 1\,,
$$
where $K$ is finite. We refer to such an extension as being \emph{finite-by-$G$}.

\begin{introtheorem}[{see Theorem~\ref{thm:E}}]\label{introB}
    Let $G$ be a group, $\widetilde{G}$ a finite-by-$G$ extension, and $\frak u\in Z^2(\widetilde{G},\mathbb T)$. Assume that every  reduced twisted group \Cs-algebra over a finite-index subgroup of $G$ is selfless with respect to the canonical trace.  Then $\Csr(\widetilde{G},\frak u)$ is pure and has stable rank one.
\end{introtheorem}

Among the group extensions covered by Theorem \ref{introB}, we highlight the following. Recall that a group is said to be acylindrically hyperbolic if it admits a non-elementary acylindrical action on a Gromov hyperbolic metric space \cite{Os16}.

\begin{introcorollary}[{see Corollary~\ref{cor:F}}]\label{MainTh}
    Let $G$ be a group, and $\frak u\in Z^2(G,\mathbb T)$. Assume that $G$ is either 
    \begin{enumerate}[leftmargin=*,label=\textup{(\roman*)}]
        \item acylindrically hyperbolic, or
        \item a lattice in ${\rm SL}(n,\mathbb R)$, for $n\geq 2$.
    \end{enumerate}
    Then the reduced twisted group \Cs-algebra $\Csr(G,\frak u)$ is pure and has stable rank one. 
\end{introcorollary}

We observe that Corollary \ref{MainTh}, part (i), removes the rapid decay assumption from Raum--Thiel--Vilalta's Theorem B \cite{RTV25}.
Moreover, note that for $n\geq 3$, ${\rm SL}(n,\mathbb Z)$ does not have the rapid decay property (see \cite[p.~55]{Ch17}). To the best of the authors' knowledge, we also establish for the first time that $\Csr({\rm SL}(n,\mathbb Z))$ has stable rank one.



\section{Actions with the $\mathrm{PHP}$ property}

\begin{definition}
Let $G$ be a discrete group. A twisted action of $G$ on a unital \Cs-algebra $A$ is a pair $(\alpha, \frak{u})$ consisting of maps $\alpha \colon G \to \mathrm{Aut}(A)$ and $\frak{u} \colon G \times G \to \mathcal{U}(A)$ satisfying
\[
\alpha_1 = \operatorname{id}_A, \qquad \alpha_s \circ \alpha_t = \operatorname{Ad}(\frak{u}_{s,t}) \circ \alpha_{st}
\]
and
\[
\frak{u}_{s,1} = \frak{u}_{1,s} = 1, \qquad \alpha_r(\frak{u}_{s,t}) \frak{u}_{r,st} = \frak{u}_{r,s} \frak{u}_{rs,t}
\]
for all $r, s, t \in G$.

The map $\frak{u}$ is referred to as a $2$-cocycle.
\end{definition}

For a twisted action as in the definition above, the reduced twisted crossed product $A \rtimes_{r,\alpha,\frak{u}} G$ is generated by an isomorphic copy of $A$ and unitaries $\lambda_g^{\frak{u}}$ satisfying 
\[
\lambda_g^{\frak{u}} a (\lambda_g^{\frak{u}})^* = \alpha_g(a), \qquad \lambda_g^{\frak{u}}\lambda_h^{\frak{u}}= \frak{u}_{g,h}\lambda_{gh}^{\frak{u}}
\]
for all $a\in A$ and $g,h\in G$.
Abusing notation, we let the \Cs-subalgebra of $A \rtimes_{r,\alpha,\frak{u}} G$ generated by the $\lambda_g^{\frak{u}}$'s be denoted by $\Cs_r(G,\frak{u})$. (Note that unless the $2$-cocycle is scalar-valued, $\Cs_r(G,\frak{u})$ is not a reduced twisted group \Cs-algebra in the usual sense.)
There is a canonical faithful conditional expectation $E\colon A \rtimes_{r,\alpha,\frak{u}} G \to A$ characterised by 
\[
E(a\lambda_g^{\frak{u}})=\begin{cases}
    a &\text{if }g=e \\
    0 &\text{otherwise.}
\end{cases}
\]

\begin{definition}\label{def:PHP}
Let $G$ be a discrete group, $A$ a unital \Cs-algebra, and $\phi$ a state on $A$. We say that a map $\alpha\colon G\to\mathrm{Aut}(A)$ has the $\mathrm{PHP}$ property (with respect to $\phi$) if for every finite $F_A \subseteq A$, finite $F_G \subseteq G$, $N\in\mathbb N$, and $\varepsilon > 0$, there exist $n\geq N$, $t_i \in G$ and $C_i \subseteq D_i \subseteq G$ for
$1 \le i \le n$ satisfying the following conditions 
\begin{enumerate}[leftmargin=*,label=\textup{(\roman*)}]
\item the subsets  
\[
\{gC_i : g \in F_G, \, i = 1, \dots, n\} \cup \{h t_j^{-1} (G \setminus D_j) : h \in F_G, \, j = 1, \dots, n\}
\]
are pairwise disjoint,
\item we have the upper bound 
\[
\sup_{g \in G} |\{i : g \in D_i\} \cup \{j : g \in t_j^{-1} (G \setminus C_j)\}| \le \varepsilon n^{1/2},
\]

\item the following averaging condition holds
\[
\max_{a\in F_A}\bigg\| \frac{1}{2n} \sum_{i=1}^n \bigl(\alpha_{t_i}(a) + \alpha_{t_i}^{-1}(a)\bigr) - \phi(a)1_A \bigg\| < \varepsilon.
\]
\end{enumerate}

If $(\alpha,\frak{u})\colon G\curvearrowright A$ is a twisted action, and $\alpha$ has the $\mathrm{PHP}$ property (with respect to a state $\phi$ on $A$), we say that $(\alpha,\frak{u})$ has the $\mathrm{PHP}$ property (with respect to $\phi$).

Moreover, if the elements $(t_i)_{i=1}^n$ can always be chosen to belong to a subgroup $H\le G$, we say that $\alpha$ (respectively, $(\alpha,\frak{u})$) has the $\mathrm{PHP}$ property relative to $H$ (with respect to $\phi$).
\end{definition}

We stress that the $\mathrm{PHP}$ property of a twisted action $(\alpha,\frak{u})$ does not directly involve explicit conditions on the cocycle $\frak u$. In particular, if $\alpha$ is a $\mathrm{PHP}$ action of $G$ on a unital \Cs-algebra $A$, then any twisted action $(\alpha,\frak u)$ automatically has property $\mathrm{PHP}$, as long as $\frak u$ is valued in the centre of $A$.

\begin{remark}
Observe that if $G$ has property $\mathrm{PHP}$ in the sense of Ozawa (see~\cite{Oz25}), then the trivial action $G\curvearrowright\mathbb C$ has the $\mathrm{PHP}$ property, as the last condition becomes vacuous. Moreover, if $A$ admits an $\alpha$-invariant state, then  $\phi$ is the unique $\alpha$-invariant state, and if $\phi$ is not $\alpha$-invariant, then there are no $\alpha$-invariant states. Note that when $A=C(X)$ with $|X|>2$ and $\alpha$ is an extremely proximal action of $G$, every $\phi\in S(A)$ fails to be $\alpha$-invariant.
\end{remark}

The proof of the following theorem is close in spirit to Ozawa's. However, realising the state on the coefficient algebra $A$, when $A\not\cong\mathbb C$, requires a new and nontrivial step in the proof.

\begin{theorem}\label{thm:crossed-product-inclusion}
Let $G$ be a discrete group, $A$ a unital \Cs-algebra with a state $\phi$,  $(\alpha,\frak{u})\colon G\curvearrowright A$ a twisted action with the $\mathrm{PHP}$ property with respect to $\phi$, and assume that $\phi \circ E$ and its restriction to $\Cs_r(G,\frak{u})$ are GNS-faithful. Then the inclusion 
\[
\Cs_r(G,\frak{u}) \subseteq (A \rtimes_{r,\alpha,\frak{u}} G,\phi \circ E)
\]
is selfless, and every intermediate \Cs-probability space is completely selfless.

Moreover, if $(\alpha,\frak{u})$ has the  $\mathrm{PHP}$ relative to a subgroup $H\le G$, and the restriction of $\phi \circ E$ to $\Cs_r(H,\frak{u}|_{H})$ is GNS-faithful, the same conclusion holds for the inclusion
\[
\Cs_r(H,\frak{u}|_{H}) \subseteq (A \rtimes_{r,\alpha,\frak{u}} G,\phi \circ E).
\]
\end{theorem}

\begin{proof}
Let $A \rtimes_{r,\alpha,\frak{u}} G$ be represented on $H_{\phi}\otimes\ell^2(G)$ in the canonical way by twisting $\pi_\phi$ with $\alpha$, where $(H_\phi,\pi_\phi,\xi_\phi)$ is the GNS representation of $\phi$, and note that the diagonal copy of $A$ under this representation commutes with $1_{H_\phi} \otimes \ell^{\infty}(G)\subseteq \mathcal B(H_{\phi}\otimes\ell^2(G))$. Moreover, we have that the vector $\xi_\phi\otimes\delta_e$ induces $\phi \circ E$ on $A \rtimes_{r,\alpha,\frak{u}} G$ and is cyclic. Hence, by uniqueness of the GNS construction we have that this representation is unitarily equivalent to the GNS representation of $\phi \circ E$. Since the latter is assumed to be faithful, we have that $A \rtimes_{r,\alpha,\frak{u}} G$ is faithfully represented on $H_{\phi}\otimes\ell^2(G)$.

For any subset $S\subseteq G$, denote by $P_S\in 1_{H_\phi} \otimes \ell^{\infty}(G)$ the projection that acts trivially on $H_\phi$ and as the indicator function of $S$ on $\ell^2(G)$. Observe that $A$ commutes with every projection of the form $P_S$, and 
\[
\lambda_g^{\frak{u}}P_S(\lambda_g^{\frak{u}})^*=P_{gS}
\]
for all $g\in G$.

Fix $\varepsilon>0$, $N\in\mathbb N$, a finite subset $F_A \subseteq A$ in the unit ball of $A$, and a finite subset $F_G\subseteq G$ containing $e\in G$. By assumption we can find $n\geq N$, $t_i\in G$, and $C_i \subseteq D_i \subseteq G$ for $1 \le i \le n$ witnessing the $\mathrm{PHP}$ property with respect to the quadruple $(\varepsilon,N,F_A,F_G)$.
Define $P_i=P_{C_i}$ and $Q_i=P_{D_i}$, 
\[
P_i' = (\lambda^{\frak{u}}_{t_i})^* Q_i^\perp \lambda^{\frak{u}}_{t_i} = P_{t_i^{-1}(G \setminus D_i)},\quad Q_i' = (\lambda^{\frak{u}}_{t_i})^* P_i^\perp \lambda^{\frak{u}}_{t_i} = P_{t_i^{-1}(G \setminus C_i)}.
\]
Since $C_i \subseteq D_i$, we have
\[
P_i' \le Q_i'.
\]
By condition (i) of the $\mathrm{PHP}$ property we have that the projections
\[
\{P_i,P_i' : 1 \le i \le n\} \subseteq 1_{H_\phi}\otimes\ell^{\infty}(G)
\]
are pairwise orthogonal.
Let
\[
T_n = \frac{1}{\sqrt{2n}} \sum_{i=1}^n \left( P_i \lambda_{t_i}^{\mathfrak{u}} + (\lambda_{t_i}^{\mathfrak{u}})^* Q_i^\perp \right) =\frac{1}{\sqrt{2n}} \sum_{i=1}^n \left( P_i \lambda_{t_i}^{\mathfrak{u}} + P_i' (\lambda_{t_i}^{\mathfrak{u}})^* \right),
\]
and 
\[
V_i=P_i\lambda_{t_i}^{\mathfrak{u}}, \qquad V_i'=P_i'(\lambda_{t_i}^{\mathfrak{u}})^* .
\]
Observe that $V_i^*V_i=(\lambda_{t_i}^{\mathfrak{u}})^*P_i\lambda_{t_i}^{\mathfrak{u}}=1-Q_i'$ and similarly $(V'_i)^*V_i'=\lambda_{t_i}^{\mathfrak{u}}P_i'(\lambda_{t_i}^{\mathfrak{u}})^*=1-Q_i$.
It follows that 
\[
T_n^* T_n = 1 - \frac{1}{2n} \sum_{i=1}^n (Q_i + Q_i'),
\]
and condition (ii) in the $\mathrm{PHP}$ property gives that 
\[
\left\| \sum_{i=1}^n (Q_i + Q_i') \right\| \le 2 \varepsilon n^{1/2},
\]
and hence
\[
\| T_n^* T_n - 1 \| \le \frac{\varepsilon}{n^{1/2}}.
\]
For an element $a\in A$ with $\Vert a\Vert \leq 1$, we have that 
\[
T_n^* a T_n = \frac{1}{2n} \sum_{i=1}^n \left( (V_i)^* a V_i + (V_i')^* a V_i' \right).
\]
Since $(\lambda_{t_i}^{\mathfrak{u}})^* a \lambda_{t_i}^{\mathfrak{u}} = \alpha_{t_i}^{-1}(a)$, it follows that
\[
(V_i)^* a V_i = \alpha_{t_i}^{-1}(a)(1 - Q_i'),\qquad (V_i')^* a V_i' = \alpha_{t_i}(a)(1 - Q_i).
\]
Thus
\[
T_n^* a T_n = \frac{1}{2n} \sum_{i=1}^n \bigl(\alpha_{t_i}(a) + \alpha_{t_i}^{-1}(a)\bigr) - \frac{1}{2n} \sum_{i=1}^n \left( \alpha_{t_i}^{-1}(a) Q_i' + \alpha_{t_i}(a) Q_i \right).
\]
By condition (ii) of the $\mathrm{PHP}$ property we have that 
\[
\left\|\frac{1}{2n} \sum_{i=1}^n \left( \alpha_{t_i}^{-1}(a) Q_i' + \alpha_{t_i}(a) Q_i \right) \right\| \le \frac{\varepsilon}{ n^{1/2}}.
\]
We now use condition (iii) of the $\mathrm{PHP}$ property to infer that, if $a\in F_A$,
\[
\|T_n^* a T_n - \phi(a) 1\| \leq \varepsilon + \frac{\varepsilon}{n^{1/2}}.
\]
Now, let $g \in F_G\setminus\{e\}$. Since $a$ commutes with the projections involved and conjugating by the unitaries $\lambda_g^{\mathfrak{u}}$ shifts the support of the projections, we can use condition (i) of the $\mathrm{PHP}$ property to show that for any $a\in A$,
\[
T_n^* a \lambda_g^{\mathfrak{u}} T_n = 0.
\]
By linearity, this passes to finite sums of elements of the form $a \lambda_g^{\mathfrak{u}}$. We can conclude that
\[
T_n^* x T_n =0
\]
for any element $x\in A \rtimes_{r,\alpha,\frak{u}} G$ with $E(x)=0$ and supported on $F_G$. Combining the previous conclusions, for any element $y\in A \rtimes_{r,\alpha,\frak{u}} G$ supported on $F_G$ and with $E(y)\in F_A$, we have that 
\[
\|T_n^* y T_n - \phi(E(y)) 1\| \leq \varepsilon + \frac{\varepsilon}{n^{1/2}}.
\]
Next, we show that $T_n+T_n^*$ is close to $\Cs_r(G,\mathfrak{u})$ in norm. First let 

\[
S_n = \frac{1}{\sqrt{2n}} \sum_{i=1}^n \left( \lambda_{t_i}^{\mathfrak{u}} + (\lambda_{t_i}^{\mathfrak{u}})^* \right) \in \Cs_r(G,\mathfrak{u}),
\]
and 
\[
R_i = Q_i - P_i, \qquad R_i' = Q_i' - P_i'.
\]
Observe that $\lambda_{t_i}^{\mathfrak{u}} R_i' (\lambda_{t_i}^{\mathfrak{u}})^* = R_i$ and a direct calculation gives that 
\[
T_n + T_n^* - S_n = -\frac{1}{\sqrt{2n}} \sum_{i=1}^n \left( R_i \lambda_{t_i}^{\mathfrak{u}} + (\lambda_{t_i}^{\mathfrak{u}})^* R_i \right).
\]
By Cauchy--Schwarz and condition (ii) in the $\mathrm{PHP}$ property, we have that 
\[
\bigg\| \sum_i R_i \lambda_{t_i}^{\mathfrak{u}} \bigg\|=\bigg\| \sum_i R_i \lambda_{t_i}^{\mathfrak{u}}R_i' \bigg\| \le \bigg\| \sum_i R_i \bigg\|^{1/2} \bigg\| \sum_i R_i' \bigg\|^{1/2} \le \varepsilon n^{1/2},
\]
and similarly for 
\[
\bigg\| \sum_i (\lambda_{t_i}^{\mathfrak{u}})^* R_i \bigg\| \leq \varepsilon n^{1/2}.
\]
Hence we conclude that 
\[
\| T_n + T_n^* - S_n \| \le \sqrt{2} \, \varepsilon.
\]
Denote by $\psi$ the vector state on $B(H_{\phi}\otimes\ell^2(G))$ defined by $\xi_\phi \otimes \delta_e$.
We want to show that $\psi(yT_nT_n^*y^*)$ is small for $y\in A \rtimes_{r,\alpha,\frak{u}} G$ supported on a finite subset of $G$. (Clearly, it is enough to show that this is true when $y$ is a contraction.)
Let $K \subseteq G$ be finite, and $P_K$ the projection defined by $K$. 
Note that, for any element $s\in G$, we have only one summand in $P_{\{s\}}T_n$ survives because the family $\{P_i,P_i'\}_{i=1}^n$ consists of pairwise orthogonal projections.
Hence,
\[
\|P_{\{s\}}T_n\|\leq\frac1{\sqrt{2n}},
\qquad
\|P_KT_n\|^2
=\|T_n^*P_KT_n\|
\leq\sum_{s\in K}\|P_{\{s\}}T_n\|^2
\leq\frac{|K|}{2n}.
\]
If $y\in A \rtimes_{r,\alpha,\frak{u}} G$ is a contraction supported on $K$, then $y^*$ is supported on $K^{-1}$ and we have that 
\[
\psi(y T_n T_n^* y^*) = \| T_n^* y^* (\xi_\phi\otimes\delta_e) \|^2 \le \| P_{K^{-1}} T_n \|^2 \| y^* (\xi_\phi\otimes\delta_e) \|^2 \le \frac{|K|}{2n}.
\]
Since the net of $T_n$'s is uniformly bounded ($T_n^*T_n\leq 1$), for a suitable cofinal ultrafilter $\mathcal U$, it yields an isometry $T \in B(H_\phi\otimes \ell^2(G))^{\mathcal{U}}$ such that
\[
T^* y T = \psi(y) 1,\quad \psi^{\mathcal{U}}(y T T^* y^*) = 0,\quad T + T^* \in \Cs_r(G,\mathfrak{u})^{\mathcal{U}}
\]
for all $y\in A \rtimes_{r,\alpha,\frak{u}} G$ with $\|E(y)\|\leq 1$, and thus for all $y$. We conclude the proof with \cite[Theorem 13]{Oz25}.

The ``moreover'' part of the statement simply follows from the fact that the extra assumption ensures that the unitaries $\lambda_{t_i}^{\frak{u}}$ belong to $\Cs_r(H,\frak{u}|_H)$, and thus $T+T^*\in \Cs_r(H,\frak{u}|_H)^{\mathcal U}$.
\end{proof}

We remark that the GNS-faithfulness assumption in the theorem above is automatic in many cases of interest.

\begin{remark}\label{rem:GNS-faithful}
Assume that $G\curvearrowright X$ is a minimal action on a compact Hausdorff space, $\alpha\colon G\curvearrowright C(X)$ the induced action, and $\frak{u}$ a $C(X)$-valued $2$-cocycle for $\alpha$.
Then, \(\rho:=\phi\circ E\) is GNS-faithful for every state \(\phi\) on \(C(X)\). 
Indeed, if an element \(z\) is in the kernel of the GNS representation of $\rho$, then
\[
0=\rho\bigl((\lambda_g^{\frak u})^*z^*z\lambda_g^{\frak u}\bigr)
=\phi\bigl(\alpha_g^{-1}(E(z^*z))\bigr)
\quad\text{for all }g\in G.
\]
If \(E(z^*z)\geq0\) were nonzero, minimality of the action and compactness of $X$ would give finitely many translates of \(E(z^*z)\) whose sum is bounded below by a positive scalar. Applying \(\phi\) contradicts previous equalities. Hence \(E(z^*z)=0\), and faithfulness of \(E\) gives \(z=0\).

Now, let \(C(Y)=\Cs_r(G,\mathfrak{u})\cap C(X)\). This is a unital $G$-invariant \Cs-subalgebra of \(C(X)\), hence an equivariant compact factor of \(X\). Moreover, $E$ restricts to a conditional expectation \(\widetilde E\colon \Cs_r(G,\mathfrak{u}) \to C(Y)\), the $2$-cocycle $\frak u$ belongs to \(C(Y)\), and a homogeneous word of degree \(g\), multiplied by \((\lambda_g^{\frak u})^*\), belongs to \(C(Y)\). Consequently
\[
\Cs_r(G,\mathfrak{u})\cong C(Y)\rtimes_{\alpha,\frak u,r}G.
\]
The $G$-action on \(Y\) is minimal, being a factor of the minimal action on \(X\), and thus the first observation above shows that $(\phi|_{C(Y)} \circ \widetilde E)$ is GNS-faithful.
Under the isomorphism $\Cs_r(G,\mathfrak{u})\cong C(Y)\rtimes_{\alpha,\frak u,r}G$ the state $(\phi|_{C(Y)} \circ \widetilde E)$ becomes $(\phi \circ E)|_{\Cs_r(G,\mathfrak{u})}$, which is therefore GNS-faithful as well.
\end{remark}


\section{Selfless inclusions from extremely proximal actions}

\begin{definition}
Let $G$ be a group, and $X$ a compact Hausdorff space with $|X|>2$. An action $G\curvearrowright X$ is extremely proximal if, for every pair of non-empty open subsets $U,V\subseteq X$, there exists $g\in G$ such that $g(X\setminus U)\subseteq V$.
\end{definition}

\begin{proposition}\label{pass}
    Let $G$ be a group, and let $G\curvearrowright X$ be a topologically free extremely proximal action. Suppose that $H\leq G$ is a finite-index subgroup. Then the restriction $H\curvearrowright X$ is a topologically free extremely proximal action.
\end{proposition}
\begin{proof}
    It is evident that $H\curvearrowright X$ remains topologically free. Since extremely proximal actions are proximal\footnote{An action $G\curvearrowright X$ is said to be \emph{proximal} if for any two points $x,y\in X$ there exists a net $(g_i)_{i \in I}\subseteq G$ such that $\lim_ig_ix=\lim_ig_iy$. If $G\curvearrowright X$ is extremely proximal and $|X|>2$, then there is a decreasing net of open neighborhoods $(U_i)_{i \in I}$ around a point $x_0\in X$ with $\bigcap_{i \in I}U_i =\{x_0\}$ (which exists because $X$ is Hausdorff) and one can find $g_i\in G$ such that $g_ix,g_iy\in U_i$ for all $i \in I$, thus ensuring proximality of the action.}, it follows from a result of Glasner, \cite[Lemma 3.2(2)]{Gla76}, that $H\curvearrowright X$ is minimal. It remains to show that $H\curvearrowright X$ is extremely proximal. For this, let $U,V\subseteq X$ be non-empty open subsets and consider a fixed right transversal $r_1,\ldots,r_m\in G$ of $H$. 
    By minimality of $H\curvearrowright X$, we may fix $x_0\in X$ and find $s_i\in H$ such that $s_i^{-1}(r_ix_0) \in V$ for $i=1,\dots,m$. In particular, this entails that $x_0\in r_i^{-1}s_iV$ for all $i$, and therefore $\bigcap_{i=1}^mr_i^{-1}s_iV$ is an open non-empty set.
    Then one can use extreme proximality of $G\curvearrowright X$ to find $g\in G$ such that $g(X\setminus \bigcap_{i=1}^mr_i^{-1}s_iV)\subseteq U$.
    As $(r_i)_{1\leq i\leq m}$ is a right transversal of $H$ in $G$, there exists an index $1\leq j\leq m$ such that $g\in Hr_j$. It follows that $gr_j^{-1}s_j\in H$ and that $gr^{-1}_js_j( X\setminus V)\subseteq U$,  because 
    $$
    g( X\setminus r^{-1}_js_jV)\subseteq g( X\setminus\cap_{i=1}^m r^{-1}_is_iV)\subseteq U.
    $$
    Hence, we conclude that $H\curvearrowright X$ is extremely proximal.
\end{proof}

The following theorem should be regarded as a generalised version of \cite[Proposition 15]{Oz25} (which, in turn, is essentially the same as \cite[Lemma 4]{delaHarpe}). However, observe that taking care of the additional condition in the $\mathrm{PHP}$ property makes the proof considerably more delicate.

\begin{theorem}\label{thm:extremely-proximal-PHP}
Let $G\curvearrowright C(X)$ be an extremely proximal, topologically free action on a compact Hausdorff space $X$ with $|X|>2$.
Then $G\curvearrowright C(X)$ has the $\mathrm{PHP}$ property with respect to any state $\phi\in S(C(X))$. Furthermore, if the restriction to the subgroup $K\leq G$ is still extremely proximal, then the action $G\curvearrowright C(X)$ has the $\mathrm{PHP}$ property with respect to any state $\phi\in S(C(X))$ and the subgroup $K$.
\end{theorem}

\begin{proof}
Fix $\varepsilon>0$ and finite subsets $F_X\subseteq C(X)$ and $F_G\subseteq G$. We may assume, without loss of generality, that $e\in F_G$ and that $F_X$ is contained in the unit ball of $C(X)$.
By the Krein--Milman Theorem, there exists $N\in\mathbb N$ and $x_i\in X$ such that 
\[
\max_{f\in F_X}\bigg| \frac{1}{N}\sum_{i=1}^Nf(x_i)-\phi(f) \bigg|<\varepsilon.
\]
Note that for any $n_0\geq N$ we can find a larger $n\geq n_0$ such that the approximation above holds true. Choose open neighbourhoods $W_i$ of $x_i$ such that 
\[
\max_{f\in F_X}|f(y)-f(x_i)|<\varepsilon
\]
for every $y\in W_i$.

Let $H=F_G^{-1}F_G$. Since the assumptions imply that $X$ has no isolated points, and the action $G\curvearrowright X$ is topologically free, we may find two points $z_i^+,z_i^- \in W_i$ for each $i=1,\dots,n$ such that the sets $\{hz_i^+,hz_i^- \mid 1\leq i\leq n,\, h\in H\}$ are all distinct. As $X$ is Hausdorff, we can choose open neighbourhoods $U_i^{\pm}$ of $z_i^{\pm}$ such that $\{hU_i^+,hU_i^- \mid 1\leq i\leq n,\, h\in H\}$ are pairwise disjoint. By extreme proximality, we can find $t_i\in G$ such that 
\[
t_i(X\setminus U_i^-)\subseteq U_i^+,\qquad t_i^{-1}(X\setminus U_i^+)\subseteq U_i^-.
\]
For a fixed $x_0\in X$, define 
\[
C_i=D_i=\{s\in G \mid sx_0\in U_i^+\}.
\]
We now verify that this choice of sets and elements satisfies the $\mathrm{PHP}$ property for $\alpha$ with respect to $\phi$.

Note that, if $r\in gC_i$ for some $g\in F_G$ and $1\leq i\leq n$, then $rx_0 \in gU_i^+$. Also, if $r\in ht_j^{-1}(G\setminus D_j)$ for some $h\in F_G$ and $1\leq j\leq n$, then $rx_0\in hU_j^-$. Since we chose a pairwise disjoint family of subsets $\{hU_i^+,hU_i^- \mid 1\leq i\leq n,\, h\in H\}$, condition (i) of the $\mathrm{PHP}$ property is satisfied.

Assume now that $s\in D_i$, or equivalently, that $sx_0\in U_i^+$, and at the same time that $s\in t_i^{-1}(G\setminus C_i)$, which means that $sx_0\in U_i^-$, then again our choice of the $U_i^{\pm}$ we have that 
\[
|\{i \mid s\in D_i\}\cup\{i \mid s\in t_i^{-1}(G\setminus C_i)\}| \leq 1.
\]
Since $\varepsilon$ is fixed, and $n$ can be arbitrarily enlarged, we can ensure that $1\leq \varepsilon n^{1/2}$, and condition (ii) of the $\mathrm{PHP}$ property is satisfied as well.

Let $y\in X$, and note that if $y\not\in U_i^-$, then $t_iy\in U_i^+$, and if $y \not\in U_i^+$, then $t_i^{-1}y \in U_i^-$.
By our choice of the sets $\{U_i^{\pm}\}_{i=1}^n$, we have that $y$ belongs at most to one of them. Hence, we have that 
\[
t_i^{-1}y\in U_i^-\subseteq W_i\quad  \text{ and } \quad  t_iy\in U_i^+\subseteq W_i
\]
for all $i$ but (at most) one index $i_0$. For all $i\neq i_0$, 
\[
\max_{f\in F_X}|f(t_iy)+f(t_i^{-1}y)-2f(x_i)| < 2\varepsilon,
\]
while for the index $i_0$, we have that 
\[
\max_{f\in F_X}|f(t_iy)+f(t_i^{-1}y)-2f(x_i)| < \varepsilon+2\max_{f\in F_X}\|f\|\leq \varepsilon+2.
\]
Then we have that 
\[
\max_{f\in F_X}\bigg| \frac{1}{2n}\sum_{i=1}^n (f(t_iy)+f(t_i^{-1}y))-\frac{1}{n}\sum_{i=1}^nf(x_i) \bigg| < \varepsilon+\frac{1}{n}.
\]
Note that, since $y\in X$ was arbitrarily picked, we have that for each $f\in F_X$, 
\begin{align*}
\bigg\| \frac{1}{2n}\sum_{i=1}^n (\alpha_{t_i^{-1}}(f)+\alpha_{t_i}(f))-\phi(f) \bigg\| &\leq \varepsilon +
\bigg\| \frac{1}{2n}\sum_{i=1}^n (\alpha_{t_i^{-1}}(f)+\alpha_{t_i}(f))-\frac{1}{n}\sum_{i=1}^nf(x_i) \bigg\| \\
&= \varepsilon +\sup_{y\in X} \bigg| \frac{1}{2n}\sum_{i=1}^n (f(t_iy)+f(t_i^{-1}y))-\frac{1}{n}\sum_{i=1}^nf(x_i) \bigg| \\
&\leq 2\varepsilon+\frac{1}{n}.
\end{align*}
Since $\varepsilon>0$ was arbitrarily small and $n$ can be made arbitrarily large, we may conclude that condition (iii) of the $\mathrm{PHP}$ property is satisfied.
\end{proof}

Combining the previous results we obtain the following corollary. 

\begin{corollary}\label{cor:extreme}
Let $X$ be a compact Hausdorff space with $|X|>2$.
For any twisted action $(\alpha,\frak{u})\colon G\curvearrowright C(X)$ such that $\alpha$ induces an extremely proximal, topologically free action on $X$, and any $\phi\in S(C(X))$, the inclusion 
\[
\Cs_r(G,\frak{u}) \subseteq (C(X)\rtimes_{r,\alpha,\frak{u}}G,\phi\circ E)
\]
is selfless, and every intermediate \Cs-probability space is completely selfless.

Moreover, if $H\le G$ is a finite index subgroup and $\frak{v}=\frak{u}|_{H}$ the restricted cocycle, the same conclusion holds for the inclusion
\[
\Cs_r(H,\frak{v}) \subseteq (C(X)\rtimes_{r,\alpha,\frak{u}}G,\phi\circ E).
\]
\end{corollary}

\begin{proof}
We only need to verify that $\phi\circ E$, $(\phi\circ E)|_{\Cs_r(G,\frak{u})}$ and $(\phi\circ E)|_{\Cs_r(H,\frak{v})}$ are GNS-faithful. However, this is automatic, as explained in Remark \ref{rem:GNS-faithful}, because the actions $G,H\curvearrowright X$ are minimal. The statement then follows from Theorem~\ref{thm:crossed-product-inclusion} and Theorem~\ref{thm:extremely-proximal-PHP}.
\end{proof}


\section{Selfless inclusions from linear groups}

We now proceed to provide yet another application of Theorem~\ref{thm:crossed-product-inclusion}. But first we need some preliminaries.

\begin{remark}
Let $n\geq 2$, and $\mathrm{PSL}(n,\mathbb R)$ act on $\mathbb P^{n-1}(\mathbb R)$ in the canonical way. Consider any Zariski-dense subgroup $\Gamma\leq \mathrm{PSL}(n,\mathbb R)$ and the restricted action $\Gamma\curvearrowright\mathbb P^{n-1}(\mathbb R)$.
We fix an $\mathbb R$-regular element $\Gamma$, which exists by \cite[Théorème A.1]{BL1993}, and let $a$ be its square. We also fix a compatible metric $d$ on the projective space. Then a lift of $a$ to $\mathrm{SL}(n,\mathbb R)$ can be chosen with real eigenvalues $\lambda_1 > \lambda_2 > \dots > \lambda_{n-1}>\lambda_n>0$.
Let $v_1,\dots,v_n$ be corresponding eigenvectors.
Set $p_+=[v_1]$ and $p_-=[v_n]$, and define the following two hyperplanes
\[
H_+ = \mathbb P(\mathrm{span}\{v_1,\dots,v_{n-1}\}),\qquad H_- = \mathbb P(\mathrm{span}\{v_2,\dots,v_n\}).
\]
Observe that, for any compact subset $K\subseteq \mathbb P^{n-1}(\mathbb R)\setminus H_-$, we have 
\[
\lim_{k\to +\infty} \max_{x\in K} \, d(a^kx,p_+)=0.
\]
In fact, if $x\in \mathbb P^{n-1}(\mathbb R)\setminus H_-$, one can write $x=[x_1v_1+\dots+x_nv_n]$ with $x_1\neq0$, and 
\[
a^kx= \bigg[v_1 + \sum_{i=2}^n \left(\frac{\lambda_i}{\lambda_1}\right)^k\frac{x_i}{x_1}v_i\bigg].
\]
By compactness of $K$, the set $\{x_i/x_1 \mid 2\leq i \leq n,\, x\in K\}$ is uniformly bounded, and by assumption $\lambda_i/\lambda_1<1$.
Hence,  we have that
\[
\lim_{k\to+\infty}\max_{x\in K} \, d(a^kx,p_+)=\lim_{k\to+\infty}\max_{x\in K} \, d\bigg( \bigg[v_1 + \sum_{i=2}^n \left(\frac{\lambda_i}{\lambda_1}\right)^k\frac{x_i}{x_1}v_i\bigg],p_+ \bigg)=0.
\]
Hence, for any open neighbourhood $U_+$ of $p_+$, and any open set $V_-$ containing $H_-$, we have that 
\[
a^k(\mathbb P^{n-1}(\mathbb R)\setminus V_-) \subseteq U_+
\]
for every sufficiently large $k\in\mathbb N$.
For the same reason, we have that for any open neighbourhood $U_-$ of $p_-$, and any open set $V_+$ containing $H_+$, 
\[
a^{-k}(\mathbb P^{n-1}(\mathbb R)\setminus V_+) \subseteq U_-
\]
for every sufficiently large $k\in\mathbb N$. 
\end{remark}

\emph{In what follows, we refer to $H_{\pm}$ and $p_{\pm}$ as in the previous remark without further reference.}

\begin{remark}\label{rem:lattices-facts}
    Let $n\geq2$, and $\Gamma \le \mathrm{PSL}(n,\mathbb R)$ be a lattice.
    Then we have the following well-known facts:
    \begin{enumerate}[leftmargin=*,label=\textup{(\roman*)}]
        \item $\Gamma$ is Zariski-dense in $\mathrm{PSL}(n,\mathbb R)$ by \cite[Theorem 3.2.5]{Zim84};
        \item the canonically induced action $\Gamma\curvearrowright\mathbb P^{n-1}(\mathbb R)$ is minimal by \cite[Lemma 8.5]{Mos73}.
    \end{enumerate}
    Moreover, it is clear that the action of $\mathrm{PSL}(n,\mathbb R)$ on $\mathbb P^{n-1}(\mathbb R)$ is transitive.
\end{remark}

\begin{lemma}\label{lem:Zariski-intersection}
Let $n\geq2$, and $\Gamma \le \mathrm{PSL}(n,\mathbb R)$ be a Zariski-dense subgroup such that the canonical action $\Gamma\curvearrowright\mathbb P^{n-1}(\mathbb R)$ is minimal.
Let $W\subseteq \mathbb P^{n-1}(\mathbb R)$ be a nonempty open subset, and $Y\subseteq \mathrm{PSL}(n,\mathbb R)$ a nonempty Zariski-open subset. 
Then, there exists $s\in \Gamma\cap Y$ such that 
\[
sp_+ \in W,\qquad sp_-\in W.
\]
\end{lemma}

\begin{proof}
Since the action $\Gamma \curvearrowright \mathbb P^{n-1}(\mathbb R)$ is minimal, we may choose $g\in \Gamma$ such that $g\cdot p_+ \in W$. 
Set $h=gag^{-1}$.
Let $h_+=g\cdot p_+$, and $H_h^-=g\cdot H_-$, and note that if $x\not\in H_h^-$ (i.e., $g^{-1}x\not\in H_-)$, then 
\[
\lim_{k\to+\infty} \, d(h^kx,h_+) = \lim_{k\to+\infty} \, d\left(ga^k(g^{-1}x)) , gp_+\right)=0,
\]
which by the same argument as before holds uniformly on compact subsets of $\mathbb P^{n-1}(\mathbb R)\setminus H_h^-$.

Let $T$ be the Zariski-closure of the cyclic group generated by $h$, and denote $Z=\mathrm{PSL}(n,\mathbb R)\setminus Y$, which is a proper Zariski-closed subset.
Let 
\[
E= \{r\in\mathrm{PSL}(n,\mathbb R) \mid T\cdot r\subseteq Z\}=\bigcap_{t\in T}t^{-1}(Z),
\]
which is Zariski-closed, and also proper because $E\subseteq Z$.
Define 
\[
R_{\pm}=\{r\in \mathrm{PSL}(n,\mathbb R) \mid rp_{\pm} \in H_h^-\},
\]
which are Zariski-closed proper subsets because the map $\mathrm{PSL}(n,\mathbb R)\ni r\mapsto rp_{\pm}$ is surjective onto $\mathbb P^{n-1}(\mathbb R)$, which follows from transitivity of $\mathrm{PSL}(n,\mathbb R)\curvearrowright\mathbb P^{n-1}(\mathbb R)$, and $H_h^-$ is a proper projective hyperplane.
Hence, we have that $E \cup R_+ \cup R_-$ is a proper Zariski-closed subset because $\mathrm{PSL}(n,\mathbb R)$ is irreducible.
Then, by Zariski-density of $\Gamma$, we may find $r\in \Gamma\setminus(E \cup R_+ \cup R_-)$.
In particular, 
\begin{equation}\label{eq:Zariski-intersection}
rp_+,\,rp_- \not\in H_h^-,\qquad (T\cdot r)\cap Y\neq \emptyset.
\end{equation}
Now, we show that for any $N\in\mathbb N$, the set $\{h^{k}r\mid k\geq N\}$ is Zariski-dense in $T\cdot r$. So fix such an $N$. 
Note that the Zariski-closure of $h\{h^kr\mid k\geq N\}=\{h^{k}r\mid k\geq N+1\}$ is contained in the Zariski-closure of $\{h^kr\mid k\geq N\}$. Since multiplying by a matrix $h$ is a Zariski-homeomorphism of $\mathrm{PSL}(n,\mathbb R)$, taking the Zariski-closure of $h\{h^kr\mid k\geq N\}$ corresponds to applying $h$ to the Zariski closure of $\{h^kr\mid k\geq N\}$.
Hence, we are left with a descending sequence of Zariski-closed sets, i.e., the Zariski-closures of $\{h^kr\mid k\geq N\} \supseteq \{h^{k}r\mid k\geq N+1\} \supseteq\dots$, which must eventually stabilise by the descending chain condition of the Zariski topology on $\mathrm{PSL}(n,\mathbb R)$.
Consequently, there exists $m_0\in\mathbb N$ such that the Zariski-closure of $\{h^{k}r\mid k\geq N+m_0\}$ is invariant under applying any integer power of $h$. But then, after applying the right integer power of $h$ to that set, we have that $h^jr$ lies in the Zariski-closure of $\{h^{k}r\mid k\geq N+m_0\}$ for all $j\in\mathbb Z$. Hence, $T \cdot r$ is contained in the Zariski-closure of $\{h^{k}r\mid k\geq N+m_0\}$, and thus any set $\{h^{k}r\mid k\geq N\}$ is Zariski-dense in $T\cdot r$.

Next, we claim that $h^kr\in Y$ for arbitrarily large $k\in\mathbb N$. 
Suppose, towards a contradiction, that there is some $N$ such that $h^kr\not\in Y$ for all $k\geq N$. Equivalently, $\{h^kr\mid k\geq N\} \subseteq Z$. Since $Z$ is Zariski-closed, the Zariski-closure of $\{h^kr\mid k\geq N\}$ remains contained in $Z$. By the observation above, we have that the Zariski-closure of $\{h^kr\mid k\geq N\}$ is Zariski-dense in $T\cdot r$, thus $T\cdot r \subseteq Z$, which contradicts \eqref{eq:Zariski-intersection}.

The properties of $h$ established at the beginning of the proof imply that
\[
\lim_{k\to+\infty} \, d\left(h^k (rp_+) , h_+\right)=0,\qquad \lim_{k\to+\infty} \, d\left(h^k (rp_-) , h_+\right)=0
\]
because $rp_{\pm}\not\in H_h^-$ (see \eqref{eq:Zariski-intersection}).
Since $h_+\in W$, for any sufficiently large $k\in\mathbb N$, $h^k (rp_+),h^k (rp_-)\in W$.
Choosing $k$ large enough so that we also have $h^kr\in Y$, we conclude that $s=h^kr$ satisfies 
\[
s\in  \Gamma\cap Y,\qquad sp_+,\,sp_-\in W.
\]
\end{proof}

The following theorem can be viewed as a generalisation of \cite[Proposition 16]{Oz25} (which follows the proof of \cite[Proposition 1]{BCdlH}).

\begin{theorem}\label{thm:PSL-PHP}
Let $n\geq2$, and $\Gamma \le \mathrm{PSL}(n,\mathbb R)$ be a Zariski-dense subgroup such that the canonical action $\Gamma\curvearrowright\mathbb P^{n-1}(\mathbb R)$ is minimal.
Then $\Gamma\curvearrowright\mathbb P^{n-1}(\mathbb R)$ has the $\mathrm{PHP}$ property with respect to any state $\phi$ on $C(\mathbb P^{n-1}(\mathbb R))$.
\end{theorem}

\begin{proof}
To simplify notation, let us denote $\mathbb P^{n-1}(\mathbb R)$ by $X$.
Fix $\varepsilon>0$, and finite sets $F_X\subseteq C(X)$, $e\in F_\Gamma\subseteq \Gamma$, and a state $\phi$ on $C(\mathbb P^{n-1}(\mathbb R))$. (We may assume, without loss of generality, that $F_X$ consists of contractions.)
We may find $N\in\mathbb N$ such that for any large enough $m\geq N$ there exist $y_1,\dots,y_m\in X$ such that 
\[
\max_{f\in F_X}\bigg| \frac{1}{m}\sum_{i=1}^mf(y_i)-\phi(f)\bigg| <\varepsilon.
\]
Find open neighbourhoods $W_i$ of $y_i$, for $i=1,\dots,m$, such that 
\[
\max_{f\in F_X}|f(x)-f(x')| < \varepsilon
\]
for every $x,x'\in W_i$.
Fix linear functionals $\varphi_{\pm}$ on $\mathbb R^n$ such that $H_{\pm}=\mathbb P(\ker\varphi_{\pm})$.
We claim that there exists $(s_i)_{i=1}^m\subseteq \Gamma$ with the following properties:
\begin{enumerate}[leftmargin=*,label=\textup{(\roman*)}]
    \item $s_ip_+,\, s_ip_-\in W_i$ for all $i=1,\dots,m$,
    \item the points $\{as_ip_{-},\, as_ip_+ \mid a\in F_\Gamma,\, i=1,\dots,m\}$ are pairwise distinct,
    \item for any subset $J$ of $\{1,\dots,m\}$,\footnote{When we write $\ker(\varphi_{\pm}\circ s_j^{-1})$, we are implicitly choosing a lift of $s_j$ to $\mathrm{SL}(n,\mathbb R)$. The choice of this lift does not matter to any of our arguments, and for this reason we abuse notation and retain the notation $s_j$.}
    \[
    \dim_{\mathbb R}\bigg(\bigcap_{j\in J}\ker(\varphi_+ \circ s_j^{-1})\bigg) = \max(n-|J|,0),
    \]
    \item for any subset $J$ of $\{1,\dots,m\}$, 
    \[
    \dim_{\mathbb R}\bigg(\bigcap_{j\in J}\ker(\varphi_- \circ s_j^{-1})\bigg) = \max(n-|J|,0),
    \]
\end{enumerate}
We proceed inductively as follows. We first note that conditions (ii)--(iv) are Zariski-open.
If we have already constructed 
\[
s_1,\dots,s_{i-1} \in\Gamma
\]
with conditions (i)--(iv), then call $Q_{i-1}=\{as_jp_{-},\, as_jp_+ \mid a\in F_\Gamma,\, j=1,\dots,i-1\}$, which consists of pairwise distinct points. Then, for any $a\in F_\Gamma$, and $q\in Q_{i-1}$, the set $B_{a,q}^{\pm}=\{s\in \mathrm{PSL}(n,\mathbb R) \mid asp_{\pm}=q\}$ is a proper Zariski-closed subset because orbit maps $\mathrm{PSL}(n,\mathbb R)\to X$ are algebraic maps and because the action is transitive. Hence, $\{s\in\mathrm{PSL}(n,\mathbb R) \mid asp_{\pm}\neq q\}$ is a Zariski-open nonempty subset of $\mathrm{PSL}(n,\mathbb R)$. If $s\in B_{a,q}^+\cap B_{b,q}^+$ for $a\neq b$ (and similarly for $B_{a,q}^-\cap B_{b,q}^-$), then $b^{-1}asp_+=sp_+$. Since the fixed point set of the nontrivial element $b^{-1}a$ is a proper Zariski-closed subset, and the action is transitive, we have that $\{s\mid asp_+\neq bsp_+\}$ is a Zariski-open nonempty subset.
Note that the image of the map $s\mapsto (sp_+,sp_-)$ is precisely $\Delta^c:=\{(x,y) \in  X\times X \mid x\neq y\}$, which is a Zariski-open set.
If $s\in B_{a,q}^+\cap B_{b,q}^-$, then $(sp_+,sp_-)$ lies in the graph of $b^{-1}a$ restricted to $\Delta^c$, which is a proper Zariski-closed subset of $\Delta^c$. Hence, $\{s\mid asp_+\neq bsp_-\}$ is a nonempty Zariski-open set. 
To ensure that conditions (iii) and (iv) are also satisfied, note that every hyperplane in $X$ can be written as $sH_+$ for some $s$ (and analogously for $H_-$).
Now, since condition (iii) holds up to the index $i-1$, 
for any index subset $J\subseteq \{1,\dots,i-1\}$ with $|J|\leq n-1$ we have that 
    \[
    \dim_{\mathbb R}\bigg(\bigcap_{j\in J}\ker(\varphi_+ \circ s_j^{-1})\bigg) = \max(n-|J|,0).
    \]
It follows that $\bigcap_{j\in J}s_j H_+$ has projective dimension $n-|J|-1$, and it is therefore nonempty. So choose a point $[w_J]\in \bigcap_{j\in J}s_j H_+$, and consider the set 
\[
B_J=\{ s\in\mathrm{PSL}(n,\mathbb R) \mid [w_J] \in sH_+ \}=\{  s\in\mathrm{PSL}(n,\mathbb R) \mid [s^{-1}w_J]\in H_+ \} .
\]
$B_J$ is the preimage of $H_+$ under the orbit map $s\mapsto [s^{-1}w_J]$.
Since the action is transitive, $B_J$ is a proper Zariski-closed subset of $\mathrm{PSL}(n,\mathbb R)$.
Hence,
\[
\{s\in\mathrm{PSL}(n,\mathbb R) \mid [w_J] \not\in sH_+ \}
\]
is Zariski-open and nonempty.
Since we have finitely many such index sets $J$, we have that 
\[
D_+=\bigcap_{J} \{s\in\mathrm{PSL}(n,\mathbb R) \mid [w_J] \not\in sH_+ \}
\]
is Zariski-open and nonempty. Moreover, for any choice of $s_i\in D_+$, and for any index set $J\subseteq \{1,\dots,i\}$ containing $i$ and with $|J|\leq n$, one has that $[w_J]\not\in s_iH_+$, and thus the projective dimension of $\bigcap_{j\in J}s_j H_+$ is exactly lowered by one when compared to $\bigcap_{j\in J\setminus\{i\}}s_j H_+$ . In other words, 
    \[
    \dim_{\mathbb R}\bigg(\bigcap_{j\in J}\ker(\varphi_+ \circ s_j^{-1}) \bigg) = \max(n-|J|,0).
    \]
One can proceed analogously to find a Zariski-open nonempty subset $D_- \subseteq \mathrm{PSL}(n,\mathbb R)$ with the property that each choice of $s_i\in D_-$ yields condition (iv) for $(s_j)_{j=1}^i$. 

Thus, intersecting all the (finitely many) Zariski-open subsets constructed above gives a new Zariski-open nonempty subset of $\mathrm{PSL}(n,\mathbb R)$ whose elements satisfy conditions (ii)--(iv) up to the index $i$; denote this set by $Y_i$. 
Then Lemma~\ref{lem:Zariski-intersection} ensures that there exists $s_i\in \Gamma\cap Y_i$ such that all conditions (i)--(iv) are satisfied simultaneously. This finishes the proof of the claim.

By our choice of the $s_i$'s, we may choose open sets $V_i^{\pm}$ containing $s_iH_{\pm}$ such that $\bigcap_{i\in I}V_i^{\pm}=\emptyset$ whenever $|I|=n$.
Hence, each $x\in X$ belongs to at most $n-1$ of the $V_i^+$ sets, and at most to $n-1$ of the $V_i^-$ sets.
By our choice of the elements $s_1,\dots,s_m \in \Gamma$, the points $\{as_ip_{\pm} \mid a\in F_\Gamma,\, i=1,\dots,m\}$ are pairwise distinct, so that we may choose mutually disjoint open neighbourhoods $O_{a,i}^{\pm}$ of $as_ip_{\pm}$.
Next, choose open neighbourhoods $U_i^{\pm}$ of $s_ip_{\pm}$, each contained in 
\[
W_i \cap V_i^{\pm} \cap \bigcap_{a\in F_\Gamma}a^{-1}O_{a,i}^{\pm}
\]
for $i=1,\dots,m$. In particular, the sets $\{aU_i^{\pm} \mid a\in F_\Gamma,\,i=1,\dots,m\}$ are mutually disjoint.
By the uniform convergence for $a^k$ (and $a^{-k}$), we may find for each $i$ a sufficiently large $k_i\in\mathbb N$ such that $t_i=s_ia^{k_i}s_i^{-1}$ satisfies 
\begin{equation}\label{eq:property-ti}
t_i(X\setminus V_i^-) \subseteq U_i^+,\qquad t_i^{-1}(X\setminus V_i^+) \subseteq U_i^-.
\end{equation}

Fix now $x_0=p_+$ and define
\[
C_i=\{g\in \Gamma\mid gx_0 \in U_i^+\} \, \subseteq \, D_i=\{g\in \Gamma\mid gx_0\in V_i^+\}.
\]
We check that these sets (and the elements $t_i$) verify condition (i) of the $\mathrm{PHP}$ property.
If $g\in aC_i$, then $a^{-1}g\in C_i$ and thus $a^{-1}gx_0\in U_i^+$. 
This implies that $gx_0\in aU_i^+$. 
On the other hand, if $g\in bt_j^{-1}(\Gamma\setminus D_j)$, then $g=bt_j^{-1}h$ for some $h\not\in D_j$. Hence, $hx_0\not\in V_j^+$. By \eqref{eq:property-ti}, $t_j^{-1}hx_0\in U_j^-$, and as a result $gx_0\in bU_j^-$. 
Since the sets $\{aU_i^{\pm} \mid a\in F_\Gamma,\, i=1,\dots,m\}$ are mutually disjoint by construction, we have verified condition (i) of the $\mathrm{PHP}$ property.

To verify the second condition, note that for any fixed $g\in \Gamma$, if $g\in D_i$, then $gx_0\in V_i^+$, and 
\[
|\{1\leq i \leq m \mid g\in D_i\}| \leq n-1.
\]
If $g\in t_i^{-1}(\Gamma\setminus C_i)$, then $t_ig\not\in C_i$, and thus $t_igx_0\not\in U_i^+$. 
Again using \eqref{eq:property-ti}, we know that $gx_0\in V_i^-$. Moreover, 
\[
|\{1\leq i\leq m \mid g\in t_i^{-1}(\Gamma\setminus C_i)\}| \leq n-1.
\]
Combining the two inequalities above, and noticing that $g$ was arbitrarily chosen, we get
\[
\sup_{g\in \Gamma}|\{1\leq i \leq m \mid g\in D_i\} \cup \{1\leq i\leq m \mid g\in t_i^{-1}(\Gamma\setminus C_i)\}|\leq 2(n-1).
\]
Hence, by taking a larger $m$, which we may, we have that 
\[
\sup_{g\in \Gamma}|\{1\leq i \leq m \mid g\in D_i\} \cup \{1\leq i\leq m \mid g\in t_i^{-1}(\Gamma\setminus C_i)\}|\leq 2(n-1) \leq \varepsilon m^{1/2},
\]
and condition (ii) of the $\mathrm{PHP}$ property is satisfied.

Finally, we check condition (iii) of the $\mathrm{PHP}$ property.
Fix $x\in X$, and note that if $x\not\in V_i^+$, then \eqref{eq:property-ti} implies that $t_i^{-1}x\in U_i^-\subseteq W_i$.
Since $s_ip_{\pm}\in W_i$ by our choice of $s_i$, then 
\[
\max_{f\in F_X}|(t_i \cdot f)(x)-f(s_ip_-)| = \max_{f\in F_X}|f(t_i^{-1}x)-f(s_ip_-)|<\varepsilon.
\]
Analogously, if $x\not\in V_i^-$, we have that $t_ix\in U_i^+\subseteq W_i$, and thus 
\[
\max_{f\in F_X}|(t_i^{-1} \cdot f)(x)-f(s_ip_+)| = \max_{f\in F_X}|f(t_ix)-f(s_ip_+)|<\varepsilon.
\]
Using that 
\[
|\{1\leq i\leq m \mid x\in V_i^{\pm}\}|\leq n-1,
\]
we know that 
\begin{align*}
\max_{f\in F_X}\bigg\| \frac{1}{2m}\sum_{i=1}^m (t_i\cdot f+{t_i}^{-1}\cdot f) - \frac{1}{2m}\sum_{i=1}^m(f(s_ip_+)+f(s_ip_-))1_{C(X)} \bigg\| \leq \varepsilon + \frac{2(n-1)}{m}.
\end{align*}
Now, since $s_ip_{\pm}\in W_i$ for all $i$, we have that 
\[
\max_{1\leq i\leq m}|f(s_ip_{\pm})-f(y_i)|<\varepsilon,
\]
and therefore 
\[
\max_{f\in F_X}\bigg| \frac{1}{2m}\sum_{i=1}^m(f(s_ip_+)+f(s_ip_-)) - \frac{1}{m}\sum_{i=1}^mf(y_i) \bigg|<\varepsilon.
\]
Since the term on the right in the equation above approximates the state $\phi$, we have that 
\[
\max_{f\in F_X}\bigg\| \frac{1}{2m}\sum_{i=1}^m (t_i\cdot f+{t_i}^{-1}\cdot f) - \phi(f)1 \bigg\| \leq 3\varepsilon + \frac{2(n-1)}{m}.
\]
Since $n$ is fixed, $m$ can be chosen to be arbitrarily large, and $\varepsilon>0$ was arbitrary, we infer that condition (iii) of the $\mathrm{PHP}$ property is satisfied as well.
\end{proof}

\begin{corollary}\label{cor:PSL}
Let $n\geq2$, and $\Gamma \le \mathrm{PSL}(n,\mathbb R)$ be a Zariski-dense subgroup such that the canonical action $\Gamma\curvearrowright\mathbb P^{n-1}(\mathbb R)$ is minimal.
Let $\frak{u}$ be any $2$-cocycle for $\alpha$, and $\phi\in S(\mathbb P^{n-1}(\mathbb R))$ any state. Then, the inclusion
\[
\Cs_r(\Gamma,\frak{u}) \subseteq (C(\mathbb P^{n-1}(\mathbb R))\rtimes_{\alpha,\frak{u},r}\Gamma,\phi \circ E)
\]
is selfless, and every intermediate \Cs-probability space is completely selfless.

In particular, the assumptions are automatic whenever $\Gamma\le \mathrm{PSL}(n,\mathbb R)$ is a lattice.
\end{corollary}

\begin{proof}
We know that $\phi\circ E$ and $(\phi\circ E)|_{\Cs_r(\Gamma,\frak{u})}$ are GNS-faithful by minimality, as explained in Remark \ref{rem:GNS-faithful}.
Then the conclusion follows from Theorem~\ref{thm:PSL-PHP} and Theorem~\ref{thm:crossed-product-inclusion}.
Any lattice in $\mathrm{PSL}(n,\mathbb R)$ satisfies the assumptions as recalled in Remark~\ref{rem:lattices-facts}.
\end{proof}


\section{Other examples of $\mathrm{PHP}$ actions}

The goal of this section is to present selfless inclusions arising from actions $\mathbb F_2\curvearrowright X$ that need not be extremely proximal. 

The terminology in the following definition is motivated by \cite[Definition 1.1]{AU22}.

\begin{definition}
    Let $G$ be a discrete group, and $\alpha\colon G\curvearrowright A$ an action on a unital \Cs-algebra. If a state $\phi\in S(A)$ satisfies that, for any finite subset $F\subseteq A$, and any $\varepsilon>0$, there exists a finitely supported probability measure $\mu\in\Prob(G)$ such that 
    \[
    \max_{a\in F}\bigg\|\sum_{g\in G}\mu(g) \alpha_g(a)-\phi(a)1_{A}\bigg\|<\varepsilon
    \]
    then we say that $\phi$ is a \emph{Powers state for $\alpha$} (or that \emph{$\alpha$ admits a Powers state $\phi$}).
\end{definition}

\begin{lemma}\label{lem:approximation-measure} Let $G$ be a discrete group, let $X$ be a compact Hausdorff space, and let $G\curvearrowright X$ be a minimal, strongly proximal action. Denote by $\alpha\colon G\curvearrowright C(X)$ the induced action. Then every state on $C(X)$ is a Powers state for $\alpha$. \end{lemma}

\begin{proof} Let $\phi$ be a state on $C(X)$, fix elements $f_{1},\ldots,f_{m}\in C(X)$ and let $F:=\{f_{1},\ldots,f_{m}\}$. Furthermore, equip $C(X)^{m}$ with the maximum norm, and let $\mathcal{C}$ be the norm closure of the convex hull of the set 
\[
\left\{ \left(\alpha_{g}(f_{1}),\ldots,\alpha_{g}(f_{m})\right)\,\mid\,g\in G\right\} .
\]

We claim that $z:=(\phi(f_{1})1_{C(X)},\ldots,\phi(f_{m})1_{C(X)})$ is contained in $\mathcal{C}$. Indeed, suppose that $z\notin\mathcal{C}$. An application of the Hahn--Banach separation theorem then implies that there exists a bounded linear functional $\eta$ on $C(X)^{m}$ with 
\begin{equation}
\text{Re}\left(\eta(z)\right)>\sup_{g\in G}\text{Re}\left(\eta\left(\alpha_{g}(f_{1}),\ldots,\alpha_{g}(f_{m})\right)\right).\label{eq:HahnBanach}
\end{equation}
For $1\leq k\leq m$ denote the embedding $C(X)$ into the $k$-th component of $C(X)^{m}$ by $\iota_{k}$. By the Riesz representation theorem, we can write $\eta\circ\iota_{k}=c_{k,1}\eta_{k,1}-c_{k,2}\eta_{k,2}+i(c_{k,3}\eta_{k,3}-c_{k,4}\eta_{k,4})$ for non-negative real numbers $c_{k,1},c_{k,2},c_{k,3},c_{k,4}\geq0$ and states $\eta_{k,1},\eta_{k,2},\eta_{k,3},\eta_{k,4}\in S(C(X))$. Consider the state
\[
\psi:=\frac{1}{4m}\sum_{k=1}^{m}\sum_{l=1}^{4}\eta_{k,l}\in S(C(X)).
\]
Since the action $G\curvearrowright X$ is minimal and strongly proximal, every point $x\in X$ admits a net $(g_{j})_{j\in J}\subseteq G$ such that $g_{j}\cdot\psi\rightarrow\delta_{x}$ in the weak$^{\ast}$-topology. By the compactness of $S(C(X))$, we may assume that the net $(g_{j}\cdot\eta_{k,l})_{j\in J}$ converges for all $1\leq k\leq m$ and $1\leq l\leq4$, in which case we have $g_{j}\cdot\eta_{k,l}\rightarrow\delta_{x}$. Hence,
\begin{eqnarray*}
& & \sup_{g\in G}\text{Re}\left(\eta\left(\alpha_{g}(f_{1}),\ldots,\alpha_{g}(f_{m})\right)\right) \geq \text{Re}\left(\eta\left(\alpha_{g_{j}^{-1}}(f_{1}),\ldots,\alpha_{g_{j}^{-1}}(f_{m})\right)\right) \\
&=& \text{Re}\left(\sum_{k=1}^{m}\left(c_{k,1}\eta_{k,1}-c_{k,2}\eta_{k,2}+i(c_{k,3}\eta_{k,3}-c_{k,4}\eta_{k,4})\right)\left(\alpha_{g_{j}^{-1}}(f_{k})\right)\right) \\
&\rightarrow& \text{Re}\left(\sum_{k=1}^{m}\left(c_{k,1}-c_{k,2}+i(c_{k,3}-c_{k,4})\right)f_{k}(x)\right) \\
&=& \text{Re}\left(\sum_{k=1}^{m}\left((\eta\circ\iota_{k})(1)\right)f_{k}(x)\right) = \text{Re}\left(\delta_{x}(h)\right),
\end{eqnarray*}
where $h:=\sum_{k=1}^{m}\left(\eta\circ\iota_{k}(1)\right)f_{k}\in C(X)$, and on the other hand
\[
\text{Re}\left(\eta(z)\right)=\text{Re}\left(\eta\left(\phi(f_{1})1_{C(X)},\ldots,\phi(f_{m})1_{C(X)}\right)\right)=\text{Re}\left(\phi(h)\right).
\]
In combination with \eqref{eq:HahnBanach}, this gives $\text{Re}\left(\phi(h)\right)>\text{Re}\left(\delta_{x}(h)\right)$ for every $x\in X$. With the Krein--Milman Theorem this leads to a contradiction. Thus, $z\in\mathcal{C}$.\\

Consequently, for every $\varepsilon>0$ there exist group elements $g_{1},\ldots,g_{n}\in G$ and non-negative real numbers $c_{1},\ldots,c_{n}\geq0$ with $\sum_{l=1}^{n}c_{l}=1$ such that 
\[
\max_{1\leq k\leq m}\left\Vert \sum_{l=1}^{n}c_{l}\alpha_{g_{l}}(f_{k})-\phi(f_{k})1_{C(X)}\right\Vert <\varepsilon.
\]
The finitely supported probability measure $\mu=\sum_{l=1}^{n}c_{l}\delta_{g_{l}}$ therefore satisfies 
\[
\max_{f\in F}\left\Vert \sum_{g\in G}\mu(g)\alpha_{g}(f)-\phi(f)1_{C(X)}\right\Vert <\varepsilon.
\]
Hence $\phi$ is a Powers state for $\alpha$. \end{proof} 

Recall that a group $G$ is acylindrically hyperbolic if it admits an acylindrical, cobounded, non-elementary action on a hyperbolic space $S$. For the relevant definitions, we refer to \cite{AbDa19,Os16}, as we will only need the space $S$ and its Gromov boundary $\partial S$ to establish the following lemma.

\begin{lemma}\label{lem:boundary-K}
    Let $G$ be an acylindrically hyperbolic group with trivial finite radical, and let $G\curvearrowright S$ be acylindrical, cobounded, and non-elementary. Then, for every nontrivial normal subgroup $K\unlhd G$, and every open nonempty disjoint subset $U,V\subseteq \partial S$, there exists $g\in K$ such that $g(\partial S\setminus V)\subseteq U$.
\end{lemma}
\begin{proof}
    We first note that $K$ still acts non-elementarily and acylindrically on $S$ \cite[Lemma 7.1]{Os16}. Furthermore, we have that the limit set $\Lambda_S(K)$ is an invariant subset of $\partial S$. Since the $G$-action on $\partial S$ is minimal \cite[Proposition 0.3]{AbDa19}, it follows that $\Lambda_S(K)=\partial S$, and the conclusion follows from minimality, north-south dynamics, and the existence of loxodromic elements in $K$.
\end{proof}

\begin{theorem}\label{thm:F2-actions}
   Let $G$ be an acylindrically hyperbolic group with trivial finite radical, and let $\alpha\colon G\curvearrowright A$ be an action on a unital \Cs-algebra that admits a Powers state $\phi\in S(A)$, such that $K:=\ker(\alpha)$ is nontrivial.
    Then $\alpha$ has the $\mathrm{PHP}$ property with respect to $\phi$.
\end{theorem}

\begin{proof}
Fix $\varepsilon>0$, $N\in\mathbb N$, a finite subset $F_A\subseteq A$, and a finite subset $E\subseteq G$. Without loss of generality, assume that $e\in E$.
First, by assumption, we may find a finitely supported measure $\mu\in\Prob(G)$ such that 
\[
    \max_{a\in F_A}\bigg\|\sum_{g\in G}\mu(g) \alpha_g(a)-\phi(a)1_{A}\bigg\|<\varepsilon/2.
\]
Enumerate the elements with nonzero weights as $r_1,\dots,r_{n_1}$.
Approximating the weights by rational numbers, followed by possibly repeating the same term multiple times, we may obtain $n_0>N$ such that  
\[
    \max_{a\in F_A}\bigg\|\frac{1}{n_0}\sum_{i=1}^{n_0} \alpha_{r_i}(a)-\phi(a)1_{A}\bigg\|<\varepsilon.
\]
Then, apply $\frac{1}{n_0}\sum_{i=1}^{n_0}\alpha_{r_i^{-1}}$ to the argument in the norm. The estimate remains true, and now the sum is indexed over a symmetric set. Therefore, we may find a new $n\ge n_0$ such that, after reindexing the elements in the sum,
\[
    \max_{a\in F_A}\bigg\|\frac{1}{2n}\sum_{i=1}^{n} (\alpha_{r_i}(a)+\alpha_{r_i^{-1}}(a))-\phi(a)1_{A}\bigg\|<\varepsilon.
\]
Fix an action $G\curvearrowright S$ as above. Note that the action induced on $\partial S$ is topologically free \cite[Proposition 0.3]{AbDa19}. Choose points $b_i^+,b_i^-\in \partial S$ for $1\le i\le n$ so that all points
\[
\{gb_i^\pm \mid g\in E,\ 1\le i\le n\}
\]  
are distinct and such that 
\[
r_ib_i^+\ne b_i^- 
\]
for all $1\le i\le n$.
This is possible by a recursive construction, where one excludes the elements of $E^{-1}E$, previously chosen translations, and the preceding $r_ib_i^+,b_i^-$ as well.
Then choose nonempty open neighbourhoods $U_i^\pm$ of $b_i^\pm$ such that 
\[
\{gU_i^\pm \mid g\in E,\ 1\le i\le n\}
\]  
are pairwise disjoint and $ r_iU_i^+\cap U_i^-=\emptyset$ for all  $1\le i\le n$. 
Since $e\in E$, the sets $U_i^\pm$ themselves are pairwise disjoint.
Lemma~\ref{lem:boundary-K} applied to $U_i^+$ and $U_i^-$ yields $u_i\in K$ such that
\[
u_i(\partial S\setminus U_i^-)\subseteq U_i^+.
\]
Set $t_i=u_i r_i u_i$, and observe that 
\begin{align*} 
t_i(\partial S\setminus U_i^-) \subseteq u_i r_i(U_i^+) \subseteq u_i(\partial S\setminus U_i^-)\subseteq U_i^+ ,
\end{align*} 
and therefore also $t_i^{-1}(\partial S\setminus U_i^+)\subseteq U_i^-$.
Now, fix any $b\in \partial S$ and define
\[ 
C_i=D_i=\{s\in  G\mid sb \in U_i^+\}.
\]
We are ready to verify the conditions of the $\mathrm{PHP}$ property.

Let $s\in gC_i$ for some $g\in E$. Then $ sb\in gU_i^+$.
If $s\in g t_i^{-1}( G\setminus D_i)$ as well, we have that 
\[
sb\in g t_i^{-1}(\partial S\setminus U_i^+) \subseteq gU_i^-.
\]  
But since the sets $\{gU_i^\pm \mid g\in E,\ 1\le i\le n\}$ are pairwise disjoint, these two conditions for $s$ cannot happen at the same time. Hence, we have condition (i) of the $\mathrm{PHP}$ property.

Fix $s\in D_i \cup t_j^{-1}(\partial S\setminus C_j)$ now. We know that either $sb\in U_i^+$ or $sb\in U_i^-$, as all the sets $U_i^\pm$ are pairwise disjoint. Hence, 
\[
\sup_{s\in  G} \left| \{i\mid s\in D_i\} \cup \{i\mid s\in t_i^{-1}( G\setminus C_i)\} \right| \le1 \le\varepsilon\sqrt n.
\]
This proves condition (ii) of the $\mathrm{PHP}$ property.

Finally, as $u_i\in K=\ker(\alpha)$, we have that 
\[
\alpha_{t_i} =\alpha_{u_i}\alpha_{r_i}\alpha_{u_i} =\alpha_{r_i}, \qquad \alpha_{t_i^{-1}}=\alpha_{r_i^{-1}}. 
\]
Thus, 
\[
\frac1{2n}\sum_{i=1}^{n} \bigl(\alpha_{t_i}+\alpha_{t_i^{-1}}\bigr)=\frac1{2n}\sum_{i=1}^{n} \bigl(\alpha_{r_i}+\alpha_{r_i^{-1}}\bigr).
\]
It follows that 
\[
\max_{a\in F_A}\left\| \frac1{2n}\sum_{i=1}^{n} \bigl(\alpha_{t_i}(a)+\alpha_{t_i^{-1}}(a)\bigr) -\varphi(f)1 \right\| <\varepsilon,
\]
which proves condition (iii) of the $\mathrm{PHP}$ property. 
\end{proof}

Combining Lemma~\ref{lem:approximation-measure} and Theorem~\ref{thm:F2-actions}, we obtain the following corollary.

\begin{corollary}\label{cor:strongly-proximal}
    Let $G$ be an acylindrically hyperbolic group whose finite radical is trivial, and $X$ a compact Hausdorff space. Let $(\alpha,\frak u)\colon G\curvearrowright C(X)$ be a twisted action whose induced action on $X$ is strongly proximal, minimal, and such that $K:=\ker(\alpha)$ is nontrivial.
    Then, for any state $\phi\in S(C(X))$, the inclusion 
    \[
    \Cs_r( G,\frak{u}) \subseteq (C(X)\rtimes_{\alpha,\frak{u},r} G, \phi \circ E)
    \]
    is selfless, and each intermediate \Cs-probability space is completely selfless. 
\end{corollary}

\begin{remark}
    We observe that the crossed products in Corollary \ref{cor:strongly-proximal}, when $|X|>1$, are traceless, for if $\phi\circ E$ is the unique trace, $\phi$ would be a $G$-invariant trace on $C(X)$, which is in contradiction with the minimal, strongly proximal action.
 It follows that the crossed products in Corollary \ref{cor:strongly-proximal} are purely infinite simple.
    
      Moreover, by \cite[Theorem 1.2]{Amrutam}, all intermediate \Cs-algebras in the inclusion 
    \[
    \Cs_r(G) \subseteq A\rtimes_{\alpha,r}G
    \]
    where $\alpha\colon G\curvearrowright A$ is as in Theorem \ref{thm:F2-actions} and $G$ is hyperbolic, must be of the form $B\rtimes_r G$ for some unital $G$-\Cs-subalgebra $B$ of $A$.
    As a result, any intermediate \Cs-probability space 
    \[
    \Cs_r(G) \subsetneq (D,(\phi\circ E)|_{D}) \subseteq (C(X)\rtimes_{\alpha,r}G,\phi \circ E),
    \]
 where $G\curvearrowright X$ is as in Corollary \ref{cor:strongly-proximal} and $G$ is hyperbolic, is purely infinite simple.
    To the best of the authors' knowledge, this has not been established before.
\end{remark}

Note that the above corollary allows us to construct many non-conjugate actions of $\mathbb F_2$ with property $\mathrm{PHP}$. 

\begin{theorem}
    There exist $2^{\aleph_0}$ strongly proximal, minimal, non-conjugate actions $\mathbb F_2\curvearrowright X_i$ such that $\mathbb F_2\curvearrowright C(X_i)$ has property $\mathrm{PHP}$.
\end{theorem}
\begin{proof}
    Note that if two actions $\phi_1:G\to {\rm Homeo}(X_1)$, $\phi_2:G\to {\rm Homeo}(X_2)$ are conjugate, then the groups $G/{\ker}(\phi_1)$ and $G/{\ker}(\phi_2)$ must be isomorphic. In order to prove the theorem, we use Ol’shanskii's construction of $2^{\aleph_0}$ nonisomorphic groups $\{G_t\}_{t\in \mathbb R}$. These groups are 2-generated Tarski monsters \cite[Theorem 28.7]{Ol91}, hence \Cs-simple by \cite[Corollary 6.6]{KaKe17}. For every $t\in \mathbb R$, fix a quotient map $\mathbb F_2\to G_t$ and let $\phi_t$ be the induced action of $\mathbb F_2$ on the Furstenberg boundary $\partial_FG_i$ of $G_i$. Then, by \Cs-simplicity, $\mathbb F_2/{\ker}(\phi_i)\cong G_i$ and the theorem follows.
\end{proof}

The following example showcases that the $\mathrm{PHP}$ can be used to produce also tracial, completely selfless crossed products.

\begin{example}\label{irrational}
    Fix an irrational number $\theta$, and let $\mathbb F_2=\langle a,b \rangle$ act on $\mathbb T$ as follows, 
\[
a\cdot z=e^{2\pi i\theta}z,\qquad b\cdot z=z\qquad \text{for all $z\in\mathbb T$}.
\]
Denote by $\alpha$ the induced action on $C(\mathbb T)$, i.e., $\alpha_a(f)(z)=f(e^{-2\pi i\theta}z)$ and $\alpha_b(f)(z)=f(z)$ for all $f\in C(\mathbb T)$ and $z\in\mathbb T$.
Let $\phi$ be the state on $C(\mathbb T)$ defined by integration against the Lebesgue measure. 
We show below that $\phi$ is a Powers state for $\alpha$. Since $\alpha$ also has nontrivial kernel, it follows from Theorem \ref{thm:F2-actions} that $C(\mathbb T)\rtimes_{\alpha,r}\mathbb F_2$ is completely selfless. Since $\phi$ happens to be $\alpha$-invariant, this crossed product is tracial.

For any monomial $f(z)=z^k$ ($0\neq k\in\mathbb Z$), and $n>0$, we have that 
\[
\bigg( \sum_{j=1}^n(\alpha_{a^j}(f) + \alpha_{a^{-j}}(f)) \bigg)(z)=
\bigg( \sum_{j=1}^n (e^{-2j\pi i k \theta } + e^{2j\pi i k \theta }) \bigg) z^k.
\]
Since the coefficient of this function is bounded from above by $\frac{2}{\vert{}1 - e^{-2\pi i k \theta}\vert{}} + \frac{2}{\vert{}1 - e^{2\pi i k \theta}\vert{}}$, independently of $n$, it follows that 
\[
\frac{1}{2n}\sum_{j=1}^n(\alpha_{a^j}(f) + \alpha_{a^{-j}}(f)) \xrightarrow{n\to\infty}0.
\]
Then, using suitable approximations by polynomials via the Stone--Weierstrass Theorem, for any finite subset $F\subseteq C(\mathbb T)$ and $\varepsilon>0$, one may ensure that there exists $n>0$ such that 
\[
\max_{f\in F}\bigg\|\frac{1}{2n}\sum_{j=1}^n(\alpha_{a^j}(f) + \alpha_{a^{-j}}(f))-\phi(f)1\bigg\| <\varepsilon.
\]
The finitely supported probability measure we are looking for is therefore given by 
\[
\mu=\frac{1}{2n}\sum_{j=1}^n(\delta_{\{a^j\}}+\delta_{\{a^{-j}\}}).
\]
\end{example}


\section{An application to finite-by-$G$ extensions}

The following lemma describes the structure of $\Csr(G,\frak u)$ when $G$ contains a finite normal subgroup. It can be proven by decomposing the algebra as a twisted crossed product (see \cite[Theorem 2.1]{Be91}) and applying the imprimitivity results of Green \cite[Theorem 2.13]{Gr80}. Such lemma has also played an important role in previous works related to strict comparison and stable rank one; see \cite{Rau25,RTV25}.  

\begin{lemma}\label{structure}
    Let $G$ be a group, and let $\widetilde{G}$ be a finite-by-$G$ extension, so that it fits in a short exact sequence of the form
    $$
    1\longrightarrow K\longrightarrow \widetilde{G}\longrightarrow G\longrightarrow 1,
    $$
    with $K$ finite. Let also $\frak u\in Z^2(\widetilde{G},\mathbb T)$. Then,  $\Csr(\widetilde{G},\frak u)$ is isomorphic to a finite direct sum of \Cs-algebras that are of the form
    $$
    \Csr(H,\nu)\otimes \mathbb M_n(\mathbb C)\otimes \mathbb M_{|G/H|}(\mathbb C),
    $$
    where $H$ is a finite-index subgroup of $G$, $n\in\mathbb N$, and $\nu\in Z^2(H,\mathbb T)$ only depends on the direct summand.
\end{lemma}

We are ready to prove Theorem~\ref{introB}.

\begin{theorem}\label{thm:E}
    Let $G$ be a group, $\widetilde{G}$ a finite-by-$G$ extension, and $\frak u\in Z^2(\widetilde{G},\mathbb T)$. Assume that every  reduced twisted group \Cs-algebra over a finite-index subgroup of $G$ is selfless with respect to the canonical trace.  Then $\Csr(\widetilde{G},\frak u)$ is pure and has stable rank one.
\end{theorem}

\begin{proof}
    By invoking Lemma \ref{structure}, we can write $\Csr(\widetilde{G},\frak u)$ as a finite direct sum of \Cs-algebras that have the form
    \begin{equation}\label{summands}
        \Csr(H,\nu)\otimes \mathbb M_n(\mathbb C)\otimes \mathbb M_{|G/H|}(\mathbb C),
    \end{equation}
    where $H$ is a finite-index subgroup of $G$, $n\in\mathbb N$, and $\nu\in Z^2(H,\mathbb T)$. By assumption, every \Cs-algebra $\Csr(H,\nu)$ is selfless, which makes each of the summands in \eqref{summands} also selfless by \cite[Theorem 4.3]{Ro25}. 
    Hence, they are simple, have a unique quasitracial state, stable rank one and strict comparison by \cite[Theorem 3.1]{Ro25}, and so are pure by \cite[Theorem 5.13]{APTV24}. The result then follows from the fact that pureness and stable rank one are preserved by direct sums.
    \end{proof}

 Then we have Corollary~\ref{MainTh} as follows.  

\begin{corollary}\label{cor:F}
    Let $G$ be a group, and $\frak u\in Z^2(G,\mathbb T)$. Assume that $G$ is either 
    \begin{enumerate}[leftmargin=*,label=\textup{(\roman*)}]
        \item acylindrically hyperbolic, or
        \item a lattice in ${\rm SL}(n,\mathbb R)$, for $n\geq 2$.
    \end{enumerate}
    Then the reduced twisted group \Cs-algebra $\Csr(G,\frak u)$ is pure and has stable rank one. 
\end{corollary}

\begin{proof}
   Assume that $G$ is acylindrically hyperbolic.
   Then $G$ is an extension of a finite normal subgroup by an acylindrically hyperbolic group with trivial finite radical (see \cite[Lemma 3.9]{MiOs15} and \cite[Lemma 1]{MiOs19}). Furthermore, it is proven in \cite{Yan25} that acylindrically hyperbolic groups with a trivial finite radical admit topologically free extreme boundary actions. Hence, the conclusion follows from Theorem \ref{introB} because groups with a topologically free extreme boundary action have property $\mathrm{PHP}$ (see \cite[Proposition 15]{Oz25} or Theorem \ref{thm:extremely-proximal-PHP}) and therefore give rise to completely selfless reduced twisted group \Cs-algebras by Theorem \ref{MainTh2}, and finite-index subgroups inherit this property thanks to Proposition \ref{pass}.
    
    For the second case, we appeal to \cite[Proposition 16]{Oz25}: Zariski-dense subgroups of ${\rm PSL}(n,\mathbb R)$ have property $\mathrm{PHP}$. Note that every finite-index subgroup of a lattice in ${\rm PSL}(n,\mathbb R)$ is also a lattice in ${\rm PSL}(n,\mathbb R)$, so it must be Zariski-dense by the Borel density theorem \cite[Theorem 3.2.5]{Zim84}. Now, denote the canonical epimorphism ${\rm SL}(n,\mathbb R)\to {\rm PSL}(n,\mathbb R)$ by $q$, and the center of ${\rm SL}(n,\mathbb R)$ by $\mathcal{Z}({\rm SL}(n,\mathbb R))$. If $G$ is any lattice in ${\rm SL}(n,\mathbb R)$, then it fits into an exact sequence of the form
    $$
    1\longrightarrow G\cap  \mathcal{Z}({\rm SL}(n,\mathbb R))\longrightarrow G\longrightarrow q(G)\longrightarrow 1.
    $$
    This shows that $G$ is of the form finite-by-$q(G)$, where $q(G)$ is a lattice in ${\rm PSL}(n,\mathbb R)$.
\end{proof}


\bigskip

\textbf{Acknowledgments.}  F.F.\ gratefully acknowledges support from the Simons Foundation Dissertation Fellowship SFI-MPS-SDF-00015100. He thanks Professor Ben Hayes for all the interesting discussions surrounding the present topic.
M.\'O.C.\ gratefully acknowledges support from the National Science Foundation under Grant No.\ NSF DMS-2144739. He thanks Professor Ben Hayes, his advisor and principal investigator of the grant, for his invaluable guidance and support.
M.P.\ was supported by the Carlsberg Foundation grant CF24-1144. 

\bigskip

\textbf{Statement on the use of AI.} ChatGPT-5.6 Sol was used for proofreading, and for finding examples of actions with property $\mathrm{PHP}$ beyond the case of extreme boundaries and projective spaces. More concretely, it provided Example \ref{irrational} and the result behind Theorem \ref{thm:F2-actions} for the case of $\mathbb F_2$. The authors then generalized the latter result to its current form.
All statements and proofs have been carefully verified and adapted by the authors, who take full responsibility for the content of the present work. 



\bibliography{bibliography}
\bibliographystyle{mybstnum}


\end{document}